\documentclass[11pt]{amsart}
\usepackage{amssymb,amsmath}
\usepackage{calrsfs}
\usepackage[mathscr]{eucal}
\usepackage{verbatim}
\def\mL{L\kern-0.08cm\char39}

\newtheorem{thm}{Theorem}
\newtheorem{cor}[thm]{Corollary}
\newtheorem{claim}{\sc Claim}[thm]

\theoremstyle{remark}

\numberwithin{thm}{section}
\begin{document}

\title{On complete metrizability of the Hausdorff metric topology}
\author{
L\'aszl\'o Zsilinszky}
\address{Department of Mathematics and Computer Science, The University of
North Carolina at Pembroke, Pembroke, NC 28372, USA}
\email{laszlo@uncp.edu}


\begin{abstract} There exists a completely metrizable bounded metrizable space $X$ with compatible metrics $d,d'$ so that the hyperspace $CL(X)$ of nonempty closed subsets of $X$ endowed with the  Hausdorff metric $H_d$, $H_{d'}$, resp.  is $\alpha$-favorable, $\beta$-favorable, resp.  in the strong Choquet game. In particular, there exists a completely metrizable bounded metric space $(X,d)$ such that $(CL(X),H_d)$ is not completely metrizable.
\end{abstract}
\maketitle
\section{Introduction} 
The {\it Hausdorff metric topology} $\tau_{H_d}$ on the hyperspace $CL(X)$ of nonempty closed subsets of a given metric space $(X,d)$ is one of the oldest and  best-studied hypertopologies due to its applicability to various areas of mathematics \cite{Au,CV,Be,KT}.  The main reason for this interest is the following well known fact \cite[\S 3.2.]{Be}:
{if $(X,d)$ is a bounded complete metric space, then $(CL(X),H_d)$ is a complete metric space}, where $H_d$ is the {\it Hausdorff metric} on $CL(X)$ defined as
\begin{equation}\label{def}
H_d(A_0,A_1)=\sup\{|d(x,A_0)-d(x,A_1)|: x\in X\}, \ \text{for} \ A_0,A_1\in CL(X),
\end{equation}
and $d(x,A)=\inf\{d(x,a): a\in A\}$ is the distance from $x\in X$ to $A\in CL(X)$.
If $d$ is  not bounded, $H_d$ is  only an infinite-valued distance, which generates the topology $\tau_{H_d}$ on $CL(X)$; moreover, since  $d'=\min\{1,d\}$ is an equivalent to $d$ bounded metric on $X$ and $\tau_{H_{d'}}=\tau_{H_d}$, we get 

\begin{thm}\label{original1}\hfill

\noindent
If $(X,d)$ is complete, then $(CL(X),\tau_{H_d})$ is completely metrizable.
\end{thm}

Various completeness-type properties of the Hausdorff metric topology are stock theorems in topology, e.g. $(CL(X),\tau_{H_d})$ is compact (resp. totally bounded) iff $X$ is \cite{Be,En}; more recently, local compactness \cite{CLP}, and cofinal completeness \cite{BDM} have been characterized for $(CL(X),\tau_{H_d})$; however, despite the above considerations and other partial results (see below), a characterization of complete metrizability of $(CL(X),\tau_{H_d})$ is unknown. 
Observe that the Hausdorff distance is sensitive to its generating metric, more precisely, $\tau_{H_{d}}=\tau_{H_{d'}}$ iff $d,d'$ are uniformly equivalent metrics on $X$ \cite[Theorem 3.3.2.]{Be}, thus, it is not automatic to argue that {\it complete metrizability} of $(X,d)$ is sufficient for complete metrizability of $(CL(X),\tau_{H_d})$ even though it is clearly necessary, since $(X,d)$ embeds as a closed subspace of $(CL(X),\tau_{H_d})$. It is the purpose of this note to demonstrate that complete metrizability of $(X,d)$, in fact, is not sufficient for complete metrizability of $(CL(X),\tau_{H_d})$, contrary to some claims in the literature \cite{BV}.

To put this question in perspective, briefly review the known results related to complete  metrizability of the Hausdorff metric topology: it was {\it Effros} \cite[Lemma]{Ef} who showed that for $(CL(X),\tau_{H_d})$ to be Polish (i.e. completely metrizable and separable), it is sufficient that $(X,d)$ is completely metrizable and totally bounded, which is in turn also necessary, since separability of $(CL(X),\tau_{H_d})$ is equivalent to total boundedness of $X$ \cite[Theorem 3.2.3.]{Be}, and $X$ sits in $(CL(X),\tau_{H_d})$ as a closed  subspace. It is possible to improve on this results using the work of {\it Costantini} \cite{Co} about another related hyperspace topology, the so-called {\it Wijsman topology}  $\tau_{W_d}$ \cite{Be}: to explain this, it is useful to view $CL(X)$ as sitting in the space $C(X)$ of real-valued continuos functions defined on $X$ via the identification $A\leftrightarrow d(\cdot,A)$, since, by (\ref{def}), $(CL(X),\tau_{H_d})$ is then a subspace of $C(X)$ with the uniform topology, while $(CL(X),\tau_{W_d})$ is a subspace of $C(X)$ with the topology of pointwise convergence. This immediately implies that $\tau_{W_d}\subseteq \tau_{H_d}$, in particular, $G_{\delta}$-subsets of  $(CL(X),\tau_{W_d})$ are $G_{\delta}$-subsets of  $(CL(X),\tau_{H_d})$ as well, which helps us to prove

\begin{thm}\hfill

\noindent
If $(X,d)$ is Polish, then $(CL(X),\tau_{H_d})$ is completely metrizable.
\end{thm}

\begin{proof} It follows from \cite{Co} that $CL(X)$ is a $G_{\delta}$-set of $(CL(\widetilde X),\tau_{W_{\widetilde d}})$, where $(\widetilde X,\widetilde d)$ is the completion of $(X,d)$. Thus, $CL(X)$  is also $G_{\delta}$ in $(CL(\widetilde X),\tau_{H_{\widetilde d}})$, therefore,
by Theorem \ref{original1}, $(CL(X),\tau_{H_d})$ is completely metrizable, since mapping $A\in CL(X)$ onto the $\widetilde X$-closure of $A$ is an isometric embedding of $(CL(X),H_d)$ into $(CL(\widetilde X),H_{\widetilde d})$ \cite{Ef}.
\end{proof}

Knowing that $\tau_{W_d}=\tau_{H_d}$ on $CL(X)$ iff $(X,d)$ is totally bounded \cite[Theorem 3.2.3.]{Be}, it is not surprising that in the above results of Effros and Costantini the Hausdorff metric and Wijsman topologies interact in studying complete metrizability of the hyperspaces, however, when a totally bounded metric is not available on $X$, i.e. when $X$ is a non-separable metric space, the two topologies have no effect on each other. Therefore the wealth of completeness results on the Wijsman topology \cite{Be1,Co1,Zs1,CP,CT} is not applicable in our case, which  demonstrates a fundamental difference between these topologies.

Since complete metrizability of a metrizable space is equivalent to its \v Cech-completeness (i.e. being $G_{\delta}$ in a compactification \cite{En}), the recent characterization of local compactness of $(CL(X),\tau_{H_d})$  by {\it Costantini, Levi, Pelant} in \cite[Corollary 15]{CLP}, as well as of the intermediary property of cofinal completeness of $(CL(X),\tau_{H_d})$ by {\it Beer, Di Maio} in \cite[Theorem 3.9.]{BDM} must be mentioned here, as they both imply complete metrizability of $(CL(X),\tau_{H_d})$.

The main results of this paper, proved in Section 3, use topological games, namely the so-called strong Choquet game and the Banach-Mazur game, which are reviewed in  Section 2, along with some relevant results about them. As mentioned in the abstract and introduction, our results will demonstrate that complete metrizability of $(X,d)$ does not guarantee the same for the Hausdorff metric topology, more specifically,  $(CL(X),\tau_{H_d})$ may not have any closed-hereditary completeness property, since it contains a closed copy of the rationals; however, we will show  this hyperspace still contains a dense completely metrizable subspace, and thus, is a Baire space.
\section{Preliminaries}

Given a metric space $(X,d)$, $A\in C(X)$ and  $\varepsilon>0$, denote by
\[B_d(A,\varepsilon)=\{x\in X: d(x,A)<\varepsilon\}
\]
the open $\varepsilon$-hull of $A$, and use $B_d(x,\varepsilon)$ instead of $B_d(\{x\},\varepsilon)$ for the open $\varepsilon$-ball about $x$. In addition to (\ref{def}), there is an equivalent definition for the  Hausdorff distance $H_d$:
\[
H_d(A_0,A_1)=\inf\{\varepsilon>0: A_0\subseteq B_d(A_1,\varepsilon) \ \text{and} \ A_1\subseteq B_d(A_0,\varepsilon)\},
\]
whenever $A_0,A_1\in CL(X)$ \cite{Be,En}. 

In the {\it strong Choquet game }  $Ch(Z)$ (cf. \cite{Ch,Ke}) players
$\alpha$ and $\beta$ take turns in choosing objects in the topological
space $Z$ with an open base $\mathcal  B$: $\beta$ starts by picking $(z_0,V_0)$
from
$\mathcal E=\{ (z,V)\in Z\times\mathcal  B: \ z\in V\}$
and $\alpha$ responds by
$U_0\in \mathcal  B$ with $z_0\in U_0\subseteq V_0$. The next choice of
$\beta$ is  $(z_1,V_1)\in \mathcal  E$ with $V_1\subseteq U_0$
and again $\alpha$ picks $U_1$ with $z_1\in U_1\subseteq V_1$ etc. Player
$\alpha$ wins the run $(z_0,V_0),U_0,\dots,(z_n,V_n),U_n,\dots$ provided
$\bigcap_{n} U_n=\bigcap_{n} V_n\neq \emptyset$;
otherwise, $\beta$ wins. 
A {\it strategy} in $Ch(Z)$ for $\alpha$ (resp. $\beta$) is a function $\sigma: \mathcal E^{<\omega}\to\mathcal B$
(resp. $\sigma: \mathcal B^{<\omega}\to\mathcal E$) such that 
\begin{gather*}
z_n\in \sigma((z_0,V_0),\dots,(z_n,V_n))\subseteq V_n \ \text{ for all} \
((z_0,V_0),\dots,(z_n,V_n))\in \mathcal E^{<\omega}\\
\text{(resp.} \ \sigma(\emptyset)=(z_0,V_0) \ \text{and} \  V_n\subseteq U_{n-1}, \text{where}\
\sigma(U_{0},\dots,U_{n-1})=(z_n,V_n)\\ \text{ for all} \
(U_{0},\dots,U_{n-1})\in \mathcal B^n, \ n\ge 1\text{).}
\end{gather*} 
A strategy $\sigma$ for $\alpha$ (resp. $\beta$) is a {\it winning strategy},
if $\alpha$ (resp. $\beta$) wins every run of $Ch(Z)$ 
compatible with $\sigma$, i.e. such that  $\sigma(z_0,V_0),\dots,(z_n,V_n))=U_n$ for all $n<\omega$
(resp. $\sigma(\emptyset)=(z_0,V_0)$ and $\sigma(U_{0},\dots,U_{n-1})=(z_n,V_n)$ for all $n\ge 1$).  
 The strong Choquet game  {\it $Ch(Z)$ is $\alpha$-, $\beta$-favorable}, respectively, provided $\alpha$, resp. $\beta$ has a winning strategy in $Ch(Z)$. This game has been studied in general topological spaces \cite{Ma,BLR,DM,De,Zs}, however, the
two fundamental results about it concern metrizable ones:
\begin{itemize}
\item {\bf Choquet} \ \cite{Ch,Ke} A metrizable space $X$ is completely metrizable if and only if $Ch(X)$ is $\alpha$-favorable.
\item  {\bf Debs-Porada-Telg\'arsky} \ \cite{De,Po,Te} A metrizable space $X$ contains a closed copy of the rationals  if and only if $Ch(X)$ is  $\beta$-favorable.
\end{itemize}

The {\it Banach-Mazur game} $BM(Z)$ (see \cite{HMC}, also referred to as  the Choquet game \cite{Ke})
is played as the strong Choquet game, except
 $\beta$'s  choice is only a nonempty open set contained in the previous
choice of $\alpha$. The notions of  $\alpha$-, $\beta$-favorability of $BM(Z)$ are defined analogously to those of $Ch(Z)$. Two key results about the Banach-Mazur game are as follows:
\begin{itemize}
\item {\bf Oxtoby} \ \cite{Ox,Wh} A metrizable space $X$ contains a dense completely metrizable subspace if and only if $BM(X)$ is $\alpha$-favorable.
\item  {\bf Oxtoby-Krom} \ \cite{Ox,HMC,Ke} A topological space $X$ is a Baire space (i.e. countable intersctions of dense open subsets of $X$ are dense)  if and only if $BM(X)$ is not $\beta$-favorable.
\end{itemize}

\section{Main results}
Our main result is as follows:

\vskip 6pt
\begin{thm}
There exists a bounded metric space $(X,d)$ such that 
\begin{enumerate}
\item $X$ is completely metrizable,
\item $(CL(X),H_d)$ contains a closed copy of the rationals; in particular, 

\noindent
$(CL(X),H_d)$ is not completely metrizable,
\item $(CL(X),H_d)$ is $\alpha$-favorable in the Banach-Mazur game; in particular, 
\noindent
$(CL(X),H_d)$ is a Baire space.
\end{enumerate}
\end{thm}

\begin{proof} (1) Consider the product space $\mathbb R^{\omega}$, where 
$\mathbb R$ has the discrete topology.
This topology is metrizable by the Baire metric
\[
d(f,g)=\frac{1}{\min \{n+1: \ f(n)\neq g(n)\}}
\]
for $f,g\in \mathbb R^{\omega}$. Denote 
$F=\{x\in\mathbb R^{\omega}: x(0)\ne 0 \ \text{and} \ x(k)=0 \ \text{for all} \ k>0\}$,
and put $X=\mathbb R^{\omega}\setminus F$. It is clear that $F$ is closed in $\mathbb R^{\omega}$, so $X$ is an open subspace of the complete space $(\mathbb R^{\omega},d)$, and hence, $(X,d)$ is completely metrizable.

\vskip 6pt

(2) By the Debs-Porada-Telg\'arsky Theorem, we need to show that 
\linebreak
$(CL(X),H_d)$ is $\beta$-favorable in the strong Choquet game:  let $\{I_{n}^0\subset \mathbb R\setminus \{0\}: n<\omega\}$ be a sequence of pairwise disjoint closed bounded intervals, and denote by $I_0$ their union. For each $t\in  I_0$ define $x_{t}^0\in X$ via
\[
x_{t}^0(k)=\begin{cases}
t, &\text{if}  \ t\in I_n^0, \ k=0 \ \text{or} \ k>n+1,\\
0, &\text{if} \ t\in I_n^0, \ 1\le k\le n+1.
\end{cases}
\]
Define  $A_0=\{x_{t}^0:  t\in  I_0\}\in CL(X)$,  $\mathbf V_0=B_{H_d}(A_0,1)$, and let $(A_0,\mathbf V_0)$ be $\beta$'s initial step in $Ch(CL(X),H_d)$. Let $\mathbf U_0=B_{H_d}(A_0,\frac{1}{n_0})$ be $\alpha$'s response, where $1\le n_0<\omega$. 
Proceeding inductively, assume we have  defined a partial run $(A_0,\mathbf V_0),\mathbf U_0,\dots,(A_m,\mathbf V_m),\mathbf U_{m}$ of the strong Choquet game in $(CL(X),H_d)$, where 
\[\mathbf U_i=B_{H_d}\left (A_i,\frac{1}{\sum\limits_{j\le i} n_j}\right )
\] 
for some $1\le n_i<\omega$ whenever $i\le m$. Moreover,  for each $1\le i\le m$ a sequence $\{I_n^i\subset I_{n_{i-1}+1}^{i-1}: n<\omega\}$  of pairwise disjoint closed bounded intervals with union $I_i$ be chosen, as well as $x_t^i\in X$ for each $t\in I_0$ so that $x_t^i=x_t^{i-1}$ whenever $t\in I_{0}\setminus I_{i}$, and for $t\in I_{i}$
\[
x_{t}^i(k)=\begin{cases}
x_{t}^{i-1}(k), &\text{if}  \ t\in I_n^{i}, \ k\le \sum\limits_{j< i}n_j,\\
0, &\text{if} \ t\in I_n^i, \ \sum\limits_{j< i}n_j< k\le 1+n+ \sum\limits_{j< i}n_j,\\
t, &\text{if} \ t\in I_n^i, \ k>1+n+ \sum\limits_{j< i}n_j.
\end{cases}
\]
Then let $A_i=\{x_t^i: t\in I_0\}$ and $\mathbf V_i=B_{H_d}(A_i,\frac{1}{1+\sum\limits_{j< i}n_j})$. Choose a sequence of pairwise disjoint closed bounded intervals $\{I_n^{m+1}\subset I_{n_m+1}^m: n<\omega\}$ with union $I_{m+1}$, and define $x_t^{m+1}=x_t^m$ for each $t\in I_0\setminus I_{m+1}$, and for $t\in I_{m+1}$ put
\[
x_{t}^{m+1}(k)=\begin{cases}
x_{t}^{m}(k), &\text{if}  \ t\in I_n^{m+1}, \ k\le \sum\limits_{i\le m} n_i,\\
0, &\text{if} \ t\in I_n^{m+1}, \ \sum\limits_{i\le m} n_i< k\le 1+n+\sum\limits_{i\le m} n_i,\\
t, &\text{if} \ t\in I_n^{m+1}, \ k>1+n+\sum\limits_{i\le m} n_i.
\end{cases}
\]
Define $A_{m+1}=\{x_t^{m+1}: t\in I_0\}$ and 
$\mathbf V_{m+1}=B_{H_d}(A_{m+1},\frac{1}{1+\sum\limits_{i\le m} n_i})$. 

\begin{claim}\label{strategy}
$\mathbf V_{m+1}\subseteq \mathbf U_{m}$.
\end{claim}

\noindent
Indeed, if $A\in \mathbf V_{m+1}$, then $A\subseteq B_d(A_{m+1},\frac{1}{1+\sum\limits_{i\le m} n_i})$, so
for all $a\in A$ there is some $x_t^{m+1}\in A_{m+1}$ with $d(a,x_t^{m+1})<\frac{1}{1+\sum\limits_{i\le m} n_i}$, which implies
\begin{equation}
a(k)=x_t^{m+1}(k) \ \text{for all} \ k\le \sum\limits_{i\le m} n_i.\label{initial}
\end{equation}
If $t\in I_0\setminus I_{m+1}$, then 
\[
d(a,x_t^{m})=d(a,x_t^{m+1})<\frac{1}{1+\sum\limits_{i\le m} n_i}<\frac{1}{\sum\limits_{i\le m} n_i};
\]
if $t\in I_{m+1}$, then $t\in I_n^{m+1}$ for some $n<\omega$. It follows from the definition of $x_t^{m+1}$, and (\ref{initial}) that  
\[
d(a,x_t^{m})\le \frac{1}{1+\sum\limits_{i\le m} n_i}<\frac{1}{\sum\limits_{i\le m} n_i},
\]
so we have $A\subseteq B_d(A_{m},\frac{1}{\sum\limits_{i\le m} n_i})$.
A similar argument shows that 
\[A_{m+1}\subseteq B_d\left (A,\frac{1}{1+\sum\limits_{i\le m} n_i}\right ) \ \text{implies} \ 
A_m\subseteq B_d\left (A,\frac{1}{\sum\limits_{i\le m} n_i}\right ),
\]
thus, $A\in \mathbf U_{m}$. As a consequence of Claim \ref{strategy}, we have that putting
\[
\sigma_{Ch}(\emptyset)=(A_{0},\mathbf V_{0}), \ \text{and} \ 
\sigma_{Ch}(\mathbf U_0,\dots,\mathbf U_m)=(A_{m+1},\mathbf V_{m+1}) \ \text{whenever} \ m<\omega,
\]
defines a strategy for player $\beta$ in the strong Choquet game on $(CL(X),H_d)$. We will be done if we prove
\begin{claim}
$\sigma_{Ch}$ is a winning strategy for $\beta$  in $Ch(CL(X),H_d)$.
\end{claim}

\noindent
To show this, consider a run
\[(A_0,\mathbf V_0),\mathbf U_0,\dots,(A_m,\mathbf V_m),\mathbf U_{m},\dots\]
of $Ch(CL(X),H_d)$ compatible with $\sigma_{Ch}$, and assume $A\in\bigcap_{m<\omega} \mathbf V_m$.  If we choose some $t\in\bigcap_{m<\omega} I_{n_m+1}^m$, note that for every $m<\omega$,
\begin{equation}
x_t^m(k)=0 \ \text{for all} \ 0<k\le 1+\sum\limits_{i\le m} n_i.\label{kukk}
\end{equation}
Since $A\in \mathbf V_0$,  there is some $a\in A$ with $d(x_t^0,a)<1$, thus,
\begin{equation}
a(0)=x_t^0(0)=t.\label{zero}
\end{equation}
Since $a\in X$, there exists $0<k$ so that
\begin{equation}
a(k)\neq 0.\label{ne}
\end{equation} 
Choose $m<\omega$ so that $k\le 1+\sum\limits_{i\le m} n_i$.
Since  $A\in \mathbf V_m$, there exists an $x_{t'}^m\in A_m$ with $d(a,x_{t'}^m)< \frac{1}{1+\sum\limits_{i\le m} n_i}$, which implies that 
\begin{gather}
x_{t'}^m(0)=a(0), \ \text{and}\label{0}
\\
x_{t'}^m(k)=a(k).\label{k}
\end{gather}
Using  (\ref{zero}),(\ref{0}) we get
\[
t'=x_{t'}^m(0)=a(0)=x_{t}^0(0)=t,
\] 
so $t'=t$. This would yield, by (\ref{k}),(\ref{kukk}), that
\[
a(k)=x_{t'}^m(k)=x_t^m(k)=0,
\] 
which contradicts (\ref{ne}). In conclusion, we got that $\bigcap_{m<\omega} \mathbf V_m=\emptyset$,  and so  $\beta$ wins in $Ch(CL(X),H_d)$.

\vskip 6pt

(3) Let $\mathbf V_0$ be $\beta$'s initial step in $BM(CL(X),H_d)$, where $\mathbf V_0=B_{H_d}(A_0,\frac{1}{n_0})$ for some $A_0\in CL(X)$ and $n_0\ge 1$. For each $a_0\in A_0$ define $x_{a_0}\in X$ via
\[
x_{a_0}(k)=\begin{cases} a_0(k), \ &\text{if} \ k<2n_0-1,\\
n_0+1, \ &\text{if} \ k\ge  2n_0-1,
\end{cases}
\]
put $C_0=\{x_{a_0}: a_0\in A_0\}$. Then $C_0\in CL(X)$, and $H_d(A_0,C_0)\le \frac{1}{2n_0}$.
Define $\mathbf U_0=B_{H_d}(C_0,\frac{1}{2n_0})$. Then $\mathbf U_0\subseteq \mathbf V_0$ (since if  $A\in \mathbf U_0$, then  $H_d(A,C_0)<\frac{1}{2n_0}$, so 
$H_d(A_0,A)\le H_d(A,C_0)+H_d(C_0,A_0)<\frac{1}{n_0}$), so  we can take $\mathbf U_0$ as $\alpha$'s first step in $BM(CL(X),H_d)$. 

Assume we have  defined a partial run $\mathbf V_0,\mathbf U_0,\dots,\mathbf V_m,\mathbf U_{m}$ of the Banach-Mazur game in $(CL(X),H_d)$, where 
\[
\mathbf V_i=B_{H_d}\left (A_i,\frac{1}{n_i}\right ) \ \text{and} \ \mathbf U_i=B_{H_d}\left (C_i,\frac{1}{2n_i}\right )
\] 
for some $2n_{i-1}\le n_i<\omega$ whenever $i\le m$ (for convenience, define $n_{-1}=\frac12$). Moreover, for each $i\le m$ let $C_i=\{x_{a_i}: a_i\in A_i\}$, where
\begin{gather}
x_{a_i}(k)=\label{popo}
\begin{cases} a_i(k), \ &\text{if} \ k< 2n_i-1,\\
1+\sum\limits_{j\le i} n_j, \ &\text{if} \ k\ge 2n_i-1.
\end{cases}
\end{gather}
Take $\mathbf V_{m+1}=B_{H_d}(A_{m+1},\frac{1}{n_{m+1}})\subseteq \mathbf U_m$.  For any
$a_{m+1}\in A_{m+1}$ define 
\[
y_{a_{m+1}}(k)=\begin{cases} a_{m+1}(k), \ &\text{if} \ k< n_{m+1},\\
2+\sum\limits_{i\le m} n_i, \ &\text{if} \ k\ge n_{m+1}.
\end{cases}
\]
Then $\{y_{a_{m+1}}: a_{m+1}\in A_{m+1}\}\in \mathbf V_{m+1}\subseteq \mathbf U_m$, so there exists $x_{a_m}\in C_m$ for some $a_m\in A_m$ so that $d(y_{a_{m+1}},x_{a_m})<\frac{1}{2n_m}$. If $n_{m+1}< 2n_m$, then 
\begin{gather*}
y_{a_{m+1}}(2n_m-1)=2+\sum\limits_{i\le m} n_i \ \text{and}\\
x_{a_m}(2n_m-1)=1+\sum\limits_{i\le m} n_i,
\end{gather*}
so $d(y_{a_{m+1}},x_{a_m})\ge \frac{1}{2n_m}$, which is impossible, thus, $n_{m+1}\ge 2n_m$. It also follows from $\mathbf V_{m+1}\subseteq \mathbf U_m$ that $H_d(A_{m+1},C_m)< \frac{1}{2n_m}$. Hence, for each $a_m\in A_m$ there exists $a_{m+1}\in A_{m+1}$ with $d(x_{a_m},a_{m+1})< \frac{1}{2n_m}$, so
\begin{equation}
a_{m+1}(k)=\label{p}
\begin{cases} a_m(k), \ &\text{if} \ k< 2n_m-1,\\
1+\sum\limits_{i\le m} n_i, \ &\text{if} \ k=2n_m-1.
\end{cases}
\end{equation}
Define  $C_{m+1}=\{x_{a_{m+1}}: a_{m+1}\in A_{m+1}\}$, where
\begin{equation}
x_{a_{m+1}}(k)=\label{x}
\begin{cases} a_{m+1}(k), \ &\text{if} \ k< 2n_{m+1}-1,\\
1+\sum\limits_{i\le m+1} n_i, \ &\text{if} \ k\ge 2n_{m+1}-1,
\end{cases}
\end{equation}
and put $\mathbf U_{m+1}=B_{H_d}(C_{m+1},\frac{1}{2n_{m+1}})$.  Note that $H_d(A_{m+1},C_{m+1})\le \frac{1}{2n_{m+1}}$, so  $\mathbf U_{m+1}\subseteq \mathbf V_{m+1}$, since if $A\in \mathbf U_{m+1}$, then  $H_d(C_{m+1},A)<\frac{1}{2n_{m+1}}$, thus,
$H_d(A_{m+1},A)\le H_d(A_{m+1},C_{m+1})+H_d(C_{m+1},A)<\frac{1}{n_{m+1}}$. This means that putting $\sigma_{BM}(\mathbf V_0,\dots,\mathbf V_m)=\mathbf U_m$ for all $m<\omega$ defines a strategy for $\alpha$ in $BM(CL(X),H_d)$.

\begin{claim}
$\sigma_{BM}$ is a winning strategy for $\alpha$  in $BM(CL(X),H_d)$.
\end{claim}

\noindent
To show this, consider a run $\mathbf V_0,\mathbf U_0,\dots,\mathbf V_{m},\mathbf U_{m},\dots$ of the Banach-Mazur game in $(CL(X),H_d)$ compatible with $\sigma_{BM}$. For any $m<\omega$ and  $a_{m}\in A_{m}$ we get an $a_{m+1}\in A_{m+1}$ satisfying (\ref{p}). Then for any $a_0\in A_0$ we can define the nonempty
\[A_1[a_0]=\{a_1\in A_1: \ a_1(k)=a_0(k) \ \text{for all} \ k<2n_0-1\}.
\]
Assume, by induction, that we have defined $A_m[a_{m-1}]\neq \emptyset$ for some $a_{m-1}\in A_{m-1}$ and $m\ge 1$. For every $a_m\in A_m[a_{m-1}]$ put
\[A_{m+1}[a_m]=\{a_{m+1}\in A_{m+1}: \ a_{m+1}(k)=a_m(k) \ \text{for all} \ k<2n_m-1\},
\]
which is nonempty by (\ref{p});  for convenience, also define $A_0[a_{-1}]=A_0$. Denote
\[
P=\{(a_m)_{m\ge 0}: a_{m}\in A_{m}[a_{m-1}] \ \text{for all} \ m\ge 0\},
\]
and for any $p=(a_m)_{m\ge 0}\in  P$ define $s_p$ as follows:  
\begin{equation}
s_p(k)=\label{z}
\begin{cases} a_{0}(k), \ &\text{if} \ k< 2n_0-1,\\
a_{m}(k), \ &\text{if} \ 2n_{m-1}-1\le k< 2n_m-1, \ m\ge 1.
\end{cases}
\end{equation}
Note, by (\ref{p}), that $s_p(2n_m-1)=1+\sum\limits_{i\le m} n_i$ for every $m<\omega$, so $s_p\in X$
for each $p\in P$.
Denote by $S$ the $X$-closure of the set $\{s_p: p\in P\}$.

Given  any $s_p\in S$, we have a sequence $p=(a_m)_{m\ge 0}\in  P$ such that
$a_i(k)=a_{i-1}(k)$ for all $1\le i\le m$ and $k<2n_i-1$, which implies by (\ref{z}) that $a_m(k)=s_p(k)$ for all $k<2n_m-1$. 

It follows that $d(s_p,A_m)\le d(s_p, a_m)\le \frac{1}{2n_m}$, so $d(s,A_m)\le \frac{1}{2n_m}<\frac{1}{n_m}$ for each $s\in S$, thus,
\begin{equation}
S\subseteq B_{H_d}\left (A_m,\frac{1}{n_m}\right ).\label{end1}
\end{equation}
Furthermore, for each $1\le i\le m$, $A_i\in \mathbf V_i\subseteq \mathbf U_{i-1}$, so for any $a_i\in A_i$ there exists $a_{i-1}\in A_{i-1}$ with $d(x_{a_{i-1}},a_i)<\frac{1}{2n_{i-1}}$, which means that $a_i(k)=x_{a_{i-1}}(k)$ for each $k\le 2n_{i-1}-1$, so by (\ref{popo}), 
\begin{equation}
a_i(k)=a_{i-1}(k) \ \text{for each} \ k< 2n_{i-1}-1;\label{end}
\end{equation}
moreover, if $i>m$ we can choose by (\ref{p}),  $a_i\in A_i$ so that (\ref{end}) is satisfied.
 It follows that  $a_i\in A_i[a_{i-1}]$ for all $1\le i$, thus, $p=(a_i)_{i\ge 0}\in P$ and $s_p(k)=a_m(k)$ for all $k<2n_m-1$. 
This implies that $d(a_m,S)\le d(a_m,s_p)\le \frac{1}{2n_m}<\frac{1}{n_m}$, so
 \begin{equation}
A_m\subseteq B_{H_d}\left (S,\frac{1}{n_m}\right ).\label{end2}
\end{equation}
In conclusion, by (\ref{end1}), (\ref{end2}) we have that  $H_d(A_m,S)<\frac{1}{n_m}$, thus,
$S\in \mathbf V_m$, which implies that $S\in \bigcap_{m<\omega} \mathbf V_{m}$, and so $\alpha$ wins.
\end{proof}

\begin{cor}
There exists a completely metrizable bounded metric space $X$ with compatible metrics $d,d'$ so that
$Ch(CL(X),H_d)$ is $\alpha$-favorable and $Ch(CL(X),H_{d'})$ is $\beta$-favorable.
\end{cor}

\bibliographystyle{amsplain}

\end{document}